\documentclass{article}
\usepackage{graphicx,subcaption}
\usepackage{mathtools}
\mathtoolsset{showonlyrefs}
\usepackage{amssymb}
\usepackage{amsmath}
\usepackage{amsthm}
\usepackage{xcolor,soul}
\usepackage{fancyhdr}
\usepackage{epstopdf}
\usepackage[colorlinks=true,citecolor=magenta,linkcolor=blue]{hyperref}
\usepackage{comment} 

\makeatletter
\let\@fnsymbol\@arabic
\makeatother

\usepackage[]{geometry}
\definecolor{greenEE}{RGB}{0, 128, 0}

\newtheorem{remark}{Remark}

\newcommand{\dd}{\mathrm{d}}
\newcommand{\norm}[1]{\left\|#1\right\|}

\newcommand{\smallO}{
  \mathchoice
    {{\scriptstyle\mathcal{O}}}% 
    {{\scriptstyle\mathcal{O}}}% \textstyle
    {{\scriptscriptstyle\mathcal{O}}}% \scriptstyle
    {\scalebox{.7}{$\scriptscriptstyle\mathcal{O}$}}%\scriptscriptstyle
  }

\usepackage{tikz}

\usetikzlibrary{positioning}
\usetikzlibrary{arrows}
\usetikzlibrary{arrows.meta}
\usetikzlibrary{calc}
\usetikzlibrary{shapes}

\graphicspath{ {./figures/} }

\title{The entry-exit function and nilpotent corrections in multidimensional slow-fast systems}
\author{Jacopo Borsotti$^1$ \\ \footnotesize{$^1$Simon Fraser University, Burnaby, Canada ({\tt jacopo\_borsotti@sfu.ca})}}
\date{}

\begin{document}

\maketitle

\begin{abstract}
    We study multidimensional entry-exit phenomena in slow-fast systems whose fast linearization contains a nilpotent component. We show that the classical spectral balance still determines the exit point as $\varepsilon \to 0$, however the nilpotent component induces an algebraic amplification invisible to the spectrum, producing a $\varepsilon|\log\varepsilon|$-shift of the exit point for $0 <\varepsilon \ll 1$. We further extend the analysis to systems with several spectral branches and show that mode dominance follows a hierarchical, lexicographic selection principle: accumulated spectral growth determines the leading order, while the nilpotent component resolves degeneracies between spectrally equivalent modes.
\end{abstract}



\section{Introduction} 

\begin{figure}[h!]
    \centering 
     \begin{tikzpicture}
 \node at (0,0) {\includegraphics[width=0.4\linewidth]{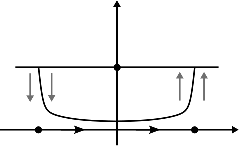}};
\node at (-0.3,2) {$x$};
\node at (-2.5,0.5) {$x=x_0$};
\node at (-2.3,-1.9) {$y_0$};
\node at (2.2,-1.9) {$y_E^\varepsilon$};
\node at (3.4,-1.9) {$y$};
  \end{tikzpicture}
    \caption{Visualization of the entry-exit map on the line $\{x=x_0\}$, $x_0$ sufficiently small, illustrating the so-called ``delayed loss of stability''. The point $(y,x)=(0,0)$ is the unique non-hyperbolic point on the critical manifold $\{x=0\}$. \label{fig:entry-exit}}
\end{figure}

The \emph{entry-exit function} \cite{MaeSch} is a standard tool of Geometric Singular Perturbation Theory (GSPT). Consider the following planar slow-fast system (the $\dot{}$ denotes the derivative with respect to $t$, whose dependency is omitted for ease of notation): 
\begin{align} \label{eq:system_entry-exit}
\begin{aligned}
    \dot{x} &= x \; f(y,\varepsilon), \\ 
    \dot{y} &= \varepsilon \; g(y,\varepsilon),
\end{aligned}
\end{align}
where $(x,y) \in \mathbb{R}^2$, $f,g$ are sufficiently regular, $g(y,0)>0$, $\textnormal{sign}(f(y,0))=\textnormal{sign}(y)$, and $0<\varepsilon \ll 1$. The entry-exit function gives, in the form of a Poincaré map, an estimate of the behavior of the orbits near the non-hyperbolic point $y=0$ of the critical manifold $\mathcal{C}_0=\{x=0\}$ \cite[Section 3.1]{Kue}. For $\varepsilon=0$, $\mathcal{C}_0$ consists of normally attracting/repelling equilibria if $y$ is negative/positive, respectively. Therefore, the intuitive explanation is to ``balance'' the contraction accumulated during the flow close to the attractive part of the critical manifold with the repulsion accumulated during the flow close to its repelling part. Refer to Figure \ref{fig:entry-exit}, it is easy to prove (see \cite[Section 5]{BorKueSen}) that an $\mathcal{O}(\varepsilon)$-approximation of the exit point $y_E^\varepsilon$ is given by the unique $y_E>0$ satisfying
\begin{equation} \label{eq:exit-point}
    \int_{y_0}^{y_E} \frac{f(y,0)}{g(y,0)}\dd y = 0.
\end{equation}
For a precise description of the orbit's behavior during the entry-exit process see \cite{Hsu}. Notice that these dynamics are mainly driven by $f(y,0)$, which corresponds to the eigenvalue describing the stability of $\mathcal{C}_0$ \cite[Definition 3.1.3]{Kue}. 

In higher dimensions, entry-exit processes become more complex. At present, no complete theory is available, however some papers studied several particular scenarios. Notably, in \cite{HsuRua} the authors considered slow-fast systems with both slow and fast variables of arbitrary finite dimension and constructed generalized entry-exit maps. However, their analysis relied on a structural assumption ensuring that the linearized fast dynamics are diagonal along the invariant manifolds, allowing contraction and expansion to be accumulated componentwise through auxiliary variables, one for each fast component. Moreover, in \cite{KakKuePopSen} the authors analyzed systems with two fast variables in which two real eigenvalues of the fast linearization intersect before contraction and expansion have balanced along either eigendirection. At the intersection, the eigenvalue has algebraic multiplicity two and geometric multiplicity one, so that the corresponding eigendirections coalesce. Specifically, the authors determine which spectral branch is effectively followed before and after the crossing. This leads to entry-exit relations in which contraction may be accumulated along one eigendirection and expansion along another, depending on the dynamics through the eigenvalue intersection. Finally, in \cite{Bor} it was shown that entry-exit phenomena arising from a particular class of three-dimensional systems can be studied focusing on just two variables, allowing for the application of \eqref{eq:exit-point}. 

In this work, we initially consider systems with $1$ slow variable and $n>1$ fast variables of the form 
\begin{align} \label{eq:system_nilpotent}
\begin{aligned}
    \dot{\mathbf{x}} &= \big(\lambda(y,\varepsilon)I_n+N\big) \mathbf{x},  \\ 
    \dot{y} &= \varepsilon \; g(y,\varepsilon),
\end{aligned}
\end{align}
where $\lambda,g$ are sufficiently regular, $g(y,0)>0$, and $\textnormal{sign}(\lambda(y,0))=\textnormal{sign}(y)$. In \eqref{eq:system_nilpotent}, $I_n$ denotes the $n \times n$ identity matrix, while $N \ne 0$ represents an $n \times n$ nilpotent matrix of index $1 <m \le n$. We will show that the exit point is given by a formula analogous to \eqref{eq:exit-point}, however with an error of order $\mathcal{O}(\varepsilon |\log \varepsilon|)$. This error, infinitely larger than the standard one in the limit $\varepsilon \to 0$, follows from the relationship between $N$ and $\mathbf{x_0}$. In Section \ref{sec:competing}, we will also analyze systems with a more general structure than \eqref{eq:system_nilpotent}. 

Importantly, the systems considered in the present work do not satisfy the structural assumptions underlying the entry-exit results discussed above. The novelty of the present work lies in isolating the effect of a persistent non-semisimple fast linearization on entry-exit dynamics. Unlike previous multidimensional extensions based on distinct fast directions or intersecting eigenvalues, the instantaneous spectrum of \eqref{eq:system_nilpotent} is unchanged by the nilpotent matrix $N$, indeed all eigenvalues related to $\{\mathbf{x}=\mathbf{0}\}$ coincide with $\lambda(y,0)$. Thus, spectral information alone predicts the same exit point as in the scalar case, while the non-semisimple structure produces an additional algebraic amplification, which increases the usual $\mathcal{O}(\varepsilon)$-error to $O(\varepsilon |\log \varepsilon|)$. This provides an explicit entry-exit correction generated purely by non-semisimplicity. Finally, we remark that errors of the order $\varepsilon |\log \varepsilon|$ were also found in \cite{HuzUld}, however they arise from degenerate planar entry-exit problems. 

\section{Main result} \label{sec:main}

Consider an initial condition $(\mathbf{x}^T, y)=(\mathbf{x_0}^T, y_0)$ with $y_0<0$ and $\norm{\mathbf{x_0}}>0$ sufficiently small. Since $y$ increases during the entry-exit process, simple calculations lead to 
\begin{equation} \label{eq:x(y)}
    \mathbf{x}(y) = \exp \left(\frac{1}{\varepsilon} \Phi(y,\varepsilon; y_0)\right) \exp \left(\frac{N}{\varepsilon} \Psi(y,\varepsilon; y_0) \right) \mathbf{x_0}, 
\end{equation}
where we defined  
\begin{equation}
    \Phi(y,\varepsilon; y_0) \coloneqq \int_{y_0}^{y} \frac{\lambda(s,\varepsilon)}{g(s,\varepsilon)} \dd s, \quad \Psi(y,\varepsilon; y_0) \coloneqq \int_{y_0}^{y} \frac{\dd s}{g(s,\varepsilon)}.
\end{equation}
Since $N$ is nilpotent with index $m$, 
\begin{equation} \label{eq:expansion}
     \exp \left(\frac{N}{\varepsilon} \Psi(y,\varepsilon; y_0) \right) \mathbf{x_0} = \mathbf{x_0} + \frac{\Psi(y,\varepsilon; y_0)}{\varepsilon} N \mathbf{x_0} + \dots + \frac{\Psi^{m-1}(y,\varepsilon; y_0)}{\varepsilon^{m-1} (m-1)!} N^{m-1} \mathbf{x_0}. 
\end{equation}
Define 
\begin{equation} \label{eq:d}
    d \coloneqq \max\{0 \le k<m \colon \; N^k \mathbf{x_0}\ne \mathbf{0}\}, 
\end{equation}
then \eqref{eq:expansion} and the regularity of $g$ imply that 
\begin{equation} \label{eq:expansion_O}
    \norm{\exp \left(\frac{N}{\varepsilon} \Psi(y,\varepsilon; y_0) \right) \mathbf{x_0}} = \frac{\Psi^{d}(y,0; y_0)}{\varepsilon^d} \frac{\norm{N^d \mathbf{x_0}}}{d!} + \mathcal{O}(\varepsilon^{1-d}); 
\end{equation}
note that \eqref{eq:expansion_O} is valid only away from $y=y_0$ since $\Psi(y_0,\varepsilon; y_0)=0$. The exit point $y_E^\varepsilon>0$ is the one satisfying $\norm{x(y_E^\varepsilon)}=\norm{\mathbf{x_0}}$. In other words, recalling \eqref{eq:x(y)}, it solves
\begin{equation} \label{eq:equality}
    \exp \left(\frac{1}{\varepsilon} \Phi(y,\varepsilon; y_0)\right) \frac{\norm{\exp \left(\frac{N}{\varepsilon} \Psi(y,\varepsilon; y_0) \right) \mathbf{x_0}}}{\norm{\mathbf{x_0}}} = 1. 
\end{equation}
Applying the logarithm to both sides of \eqref{eq:equality} and exploiting \eqref{eq:expansion_O} and the regularity of $\lambda$, we get 
\begin{equation} \label{eq:F}
    F(y,\varepsilon; y_0) \coloneqq \frac{\Phi(y,0; y_0)}{\varepsilon} + d \log \left(\frac{1}{\varepsilon}\right) + \mathcal{O}(1) = 0. 
\end{equation}
The dominant term of $F(y,\varepsilon; y_0)$ is the one depending on $\Phi$, which admits a unique $y_E>0$ such that $\Phi(y_E,0;y_0)=0$. Since the other terms are of the order at most $\mathcal{O}(|\log \varepsilon|)$, for sufficiently small $\varepsilon$, there exists a unique zero of $F(y,\varepsilon; y_0)$ in a $\mathcal{O}(\varepsilon|\log \varepsilon|)$-neighborhood of $y_E$. Such zero represents the true exit point $y_E^\varepsilon$, which is approximated by $y_E$ and coincides with it in the limit $\varepsilon \to 0$. Specifically, expanding $\Phi$ around $y_E$, it is straightforward to get 
\begin{equation} \label{eq:real_exit}
    y_E^\varepsilon = y_E - d \frac{g(y_E, 0)}{\lambda(y_E, 0)} \varepsilon \log \left(\frac{1}{\varepsilon}\right) + \mathcal{O}(\varepsilon). 
\end{equation}

\begin{remark}
    If $n=1$, then \eqref{eq:system_nilpotent} reduces to \eqref{eq:system_entry-exit} because $N=0$. In this case, $d=0$ and $y_E^\varepsilon$ is given by an $\mathcal{O}(\varepsilon)$-approximation of $y_E$. Hence, we correctly recover the standard result. Notice also that, when $N \ne 0$, we can have $d=0$ for some $\mathbf{x_0}$, however not for all initial conditions. In this case, we get again an $\mathcal{O}(\varepsilon)$-error (see Section \ref{sec:example_1}). 
\end{remark}

Eq. \eqref{eq:x(y)} isolates two contributions. Indeed, the first exponential depending on $\Phi$ is related to the spectrum (i.e., to $\lambda$), while the second one is related to the nilpotent matrix $N$. Eq. \eqref{eq:F} demonstrates that the eigenvalues play the most important role, followed by the relation between $N$ and the initial conditions (i.e., $d$ \eqref{eq:d}). 

Finally, we remark that, unlike what happens for planar systems, there may exist multiple points within a $\varepsilon$-neighborhood of $y_0$ at which $\norm{\mathbf{x}(y)}=\norm{\mathbf{x}(y_0)}$. Indeed, the asymptotic expansion \eqref{eq:expansion_O} is not valid uniformly in this region. However, these points are not relevant to the present analysis, as they arise from the initial transient dynamics rather than from the delayed loss of stability associated with the entry–exit phenomenon. 

\section{Competing spectral families and nilpotent corrections} \label{sec:competing}

Consider the following generalization of \eqref{eq:system_nilpotent}: 
\begin{align} \label{eq:system_nilpotent_general}
\begin{aligned}
    \dot{\mathbf{x}} &= \left(\bigoplus_{j=1}^r\left(\lambda_j(y,\varepsilon)I_{n_j}+N_j\right)\right) \mathbf{x},  \\ 
    \dot{y} &= \varepsilon \; g(y,\varepsilon),
\end{aligned}
\end{align}
where each $\lambda_j$ is sufficiently regular and such that $\textnormal{sign}(\lambda_j(y,0))=\textnormal{sign}(y-\mu_j)$ with $\mu_j > y_0$, $I_{n_j}$ denotes the $n_j \times n_j$ identity matrix with $\sum_{j=1}^r n_j=n$, and the $N_j$ are nilpotent $n_j \times n_j$ matrices of index $1 \le m_j \le n_j$ (such that there exists $j$ for which $N_j \ne 0$, otherwise we would have a standard entry-exit process for each variable). Let us write $\mathbf{x}^T=(\mathbf{x_1}^T, \dots, \mathbf{x_r}^T)$, with $\mathbf{x_j} \in \mathbb{R}^{n_j}$. The analysis of Section \ref{sec:main} can be applied to each $\mathbf{x_j}$ separately. In particular, 
\begin{equation} \label{eq:norm_as}
    \norm{\mathbf{x_j}(y)} \sim C_j \exp \left(\frac{1}{\varepsilon} \Phi_j(y,0; y_0)\right) \varepsilon^{-d_j},  
\end{equation}
where the $C_j>0$ depend on $y$, the initial conditions, and the system under consideration (recall that \eqref{eq:norm_as} is valid only away from $y=y_0$ since, if $d_j>0$, $C_j \to 0$ as $y \to y_0$). Since
\begin{equation}
    \norm{\mathbf{x}(y)}^2 = \sum_{j=1}^r \norm{\mathbf{x_j}(y)}^2 \sim \sum_{j=1}^r C_j^2 \exp \left(\frac{2}{\varepsilon} \Phi_j(y,0; y_0)\right) \varepsilon^{-2d_j},  
\end{equation}
$\displaystyle y_E>\min_{j \in \mathcal{J}} \mu_j$ is given by 
\begin{equation} \label{eq:exit_point_multiple}
    \max_{j \in \mathcal{J}} \Phi_j(y,0; y_0)=0, 
\end{equation}
where $\emptyset \ne \mathcal{J} = \{j \colon \; \mathbf{x_{j,0}} \ne \mathbf{0}\}$. If the dominant $\Phi_j$ corresponds to $d_j>0$, then $y_E$ represents an $\varepsilon |\log \varepsilon|$-approximation of $y_E^\varepsilon$, otherwise an $\varepsilon$-approximation. 

Consider now $i \ne j$, away from $y=y_0$, 
\begin{equation} \label{eq:ratio}
    \frac{\norm{\mathbf{x_i}(y)}}{\norm{\mathbf{x_j}(y)}} \sim \frac{C_i}{C_j} \exp \left(\frac{1}{\varepsilon} \big(\Phi_i(y,0; y_0) - \Phi_j(y,0; y_0)\big) \right) \varepsilon^{-(d_i-d_j)}, 
\end{equation}
hence we get a lexicographic order: first $\Phi_j$, then $d_j$ (similarly to what concluded in Section \ref{sec:main}). This means that as long as $\Phi_i(y,0; y_0) > \Phi_j(y,0; y_0)$, $\mathbf{x_i}$ dominates, while if $\Phi_i(y,0; y_0) < \Phi_j(y,0; y_0)$ then $\mathbf{x_j}$ dominates. At a crossing point $y_c$ such that $\Phi_i(y_c,0; y_0) = \Phi_j(y_c,0; y_0)$ one has to look at $d_i$ and $d_j$. If also the latter coincide neither $\mathbf{x_i}$ nor $\mathbf{x_j}$ dominates completely because the ratio \eqref{eq:ratio} would be of order $\mathcal{O}(1)$. However, this is precise only in the limit $\varepsilon \to 0$. The real crossing point $y_c^\varepsilon$ is the one satisfying $\Phi_i(y_c^\varepsilon,\varepsilon; y_0) = \Phi_j(y_c^\varepsilon,\varepsilon; y_0)$, which, assuming that $\lambda_i(y_c,\varepsilon) \ne \lambda_j(y_c,\varepsilon)$, can be written as $y_c^\varepsilon=y_c+\mathcal{O}(\varepsilon)$. More precisely, define the related switching point $y_{sw}^\varepsilon$ as $\norm{\mathbf{x_i}(y_{sw}^\varepsilon)}=\norm{\mathbf{x_j}(y_{sw}^\varepsilon)}$, i.e., $y_{sw}^\varepsilon$ is the point where the dominating variable changes. We have
\begin{equation}
\begin{split}
    0 = \log \left(\frac{\norm{\mathbf{x_i}(y_{sw}^\varepsilon)}}{\norm{\mathbf{x_j}(y_{sw}^\varepsilon)}}\right) &= \frac{\Phi_i(y_{sw}^\varepsilon,\varepsilon; y_0) - \Phi_j(y_{sw}^\varepsilon,\varepsilon; y_0)}{\varepsilon} + (d_i-d_j) \log \left(\frac{1}{\varepsilon}\right) + \mathcal{O}(1) \\ 
    &= \frac{\lambda_i(y_c^\varepsilon,\varepsilon) - \lambda_j(y_c^\varepsilon,\varepsilon)}{\varepsilon g(y_c^\varepsilon, \varepsilon)} (y_{sw}^\varepsilon - y_c^\varepsilon) + (d_i-d_j) \log \left(\frac{1}{\varepsilon}\right) + \mathcal{O}(1),  \\
\end{split}
\end{equation}
hence 
\begin{equation} \label{eq:sw_domination}
\begin{split}
    y_{sw}^\varepsilon &= y_c^\varepsilon - (d_i - d_j) \frac{g(y_c^\varepsilon, \varepsilon)}{\lambda_i(y_c^\varepsilon,\varepsilon) - \lambda_j(y_c^\varepsilon,\varepsilon)} \varepsilon \log \left(\frac{1}{\varepsilon}\right) + \mathcal{O}(\varepsilon) \\
    &= y_c - (d_i - d_j) \frac{g(y_c, 0)}{\lambda_i(y_c,0) - \lambda_j(y_c,0)} \varepsilon \log \left(\frac{1}{\varepsilon}\right) + \mathcal{O}(\varepsilon). 
\end{split}
\end{equation}
In conclusion, when $\Phi_i-\Phi_j \in \mathcal{O}(1)$ the block with the largest $\Phi$ dominates. On the other hand, when $\Phi_i-\Phi_j \in \mathcal{O}(\varepsilon |\log \varepsilon|)$ the dominating block depends also on $d_i$ and $d_j$. 

Finally, similarly to before, in a $\varepsilon$-neighborhood of $y_0$ multiple switching points may exist. Again, they are not relevant for our analysis. 

\section{Examples} \label{sec:examples}

\subsection{Comparison with explicit exit points} \label{sec:example_1}

Let $\mathbf{x}=(x_1, x_2)^T$ and consider 
\begin{equation} \label{eq:example_1}
    \begin{pmatrix}
        \dot{x}_1 \\ \dot{x}_2 
    \end{pmatrix}=
    \begin{pmatrix}
        y+\varepsilon & 1 \\ 0 & y+\varepsilon 
    \end{pmatrix}
    \begin{pmatrix}
        x_1 \\ x_2 
    \end{pmatrix}, \quad
    \dot{y}=\varepsilon, 
\end{equation}
and the initial conditions $y=y_0<0$, $\mathbf{x_0}^{(1)}=(1,0)^T$, and $\mathbf{x_0}^{(2)}=(0,1)^T$. 

Let us start from $\mathbf{x_0}^{(1)}=(1,0)^T$. It is easy to see that $x_2 \equiv 0$, while $x_1(y)=\exp(y-y_0+(y^2-y_0^2)/2\varepsilon)$, hence the real exit point is $y_E^\varepsilon=-y_0-2\varepsilon$. The positive zero of $\Phi(y_,0;y_0)=(y^2-y_0^2)/2$ is $y_E=-y_0$, which represents an $\mathcal{O}(\varepsilon)$-approximation of $y_E^\varepsilon$. Indeed, since $d=0$, we recover the classical error; we basically have a planar entry-exit phenomenon since $x_2 \equiv 0$. 

Consider now $\mathbf{x_0}^{(2)}=(0,1)^T$. Direct calculations show that 
\begin{equation}
    \norm{\mathbf{x}(y)}^2=\exp\left(2(y-y_0)+\frac{y^2-y_0^2}{\varepsilon}\right)\left(1+\frac{(y-y_0)^2}{\varepsilon^2}\right)
\end{equation}
and that the real exit point is 
\begin{equation}
    y_E^\varepsilon=-y_0-\frac{\varepsilon}{y_0} \log \varepsilon - 2 \varepsilon+\frac{\varepsilon}{y_0}\log (-2y_0) + \smallO(\varepsilon),  
\end{equation}
in line with \eqref{eq:real_exit}. Now the approximation given by $y_E=-y_0$ provides an $\mathcal{O}(\varepsilon|\log \varepsilon|)$-error, indeed $d=1$. 

\subsection{A numerical validation}

Let $\mathbf{x}^T=(\mathbf{x_1}^T, \mathbf{x_2}^T)$ with $\mathbf{x_1} \in \mathbb{R}^3$ and $\mathbf{x_2} \in \mathbb{R}^2$. Consider the system 
\begin{equation} \label{eq:example_2}
    \begin{pmatrix}
        \dot{\mathbf{x}}_{\mathbf{1}} \\ \dot{\mathbf{x}}_{\mathbf{2}}
    \end{pmatrix}=
    \begin{pmatrix}
        \lambda_1(y,\varepsilon)I_3+N_1 & 0 \\ 0 & \lambda_2(y,\varepsilon)I_2+N_2 
    \end{pmatrix} 
    \begin{pmatrix}
        \mathbf{x_1} \\ \mathbf{x_2} 
    \end{pmatrix}, \quad
    \dot{y}=\varepsilon(1+\varepsilon), 
\end{equation}
where 
\begin{equation}
    \lambda_1(y,\varepsilon)=y+\varepsilon, \quad \lambda_2(y,\varepsilon)=2y+\frac{1}{4}-\varepsilon
\end{equation}
and 
\begin{equation}
    N_1 = \begin{pmatrix}
        0 & 1 & 0 \\ 0 & 0 & 1 \\ 0 & 0 & 0
    \end{pmatrix}, \quad N_2=  \begin{pmatrix}
        0 & 1  \\ 0 & 0 
    \end{pmatrix}, 
\end{equation}
and the initial condition $y_0=-1$, $\mathbf{x_{1,0}}=(0,0,1)^T$, and $\mathbf{x_{2,0}}=(0,1)^T$. 

Now, 
\begin{align}
    \max \{\Phi_1(y,0; -1), \Phi_2(y,0; -1)\} = 
    \begin{cases}
        \Phi_1(y,0; -1) =  \frac{y^2-1}{2} &\textnormal{if} \; -1 \le y < \frac{1}{2}, \\ 
        \Phi_2(y,0;-1)=y^2-1+\frac{y+1}{4} \quad & \textnormal{if} \; y \ge \frac{1}{2}, 
    \end{cases} 
\end{align}
hence \eqref{eq:exit_point_multiple} predicts an $\varepsilon|\log \varepsilon|$-approximation of $y_E^\varepsilon$ given by $y_E=3/4$, specifically from \eqref{eq:real_exit} we have $y_E^\varepsilon=3/4-4/7\varepsilon|\log\varepsilon|+\mathcal{O}(\varepsilon)$. Moreover, \eqref{eq:sw_domination} implies that $y_{sw}^\varepsilon=1/2+4/3\varepsilon |\log\varepsilon|+\mathcal{O}(\varepsilon)$, for which $\norm{\mathbf{x_1}}=\norm{\mathbf{x_2}}$. Finally, in the limit $\varepsilon \to 0$,  $\mathbf{x_1}$ exponentially dominates as long as $y<1/2=y_c$, while $\mathbf{x_2}$ exponentially dominates when $y>1/2$. At $y=1/2$, since $d_1=2>1=d_2$, we have $\norm{\mathbf{x_1}}/\norm{\mathbf{x_2}} \sim K\varepsilon^{-1}$ for some $K>0$. Interestingly, as $\varepsilon \to 0$, although $y_{sw}^\varepsilon\to y_c$, at $y=y_c$ the ratio $\norm{\mathbf{x_1}}/\norm{\mathbf{x_2}}$ diverges as $\varepsilon^{-1}$. 

\begin{figure}[h!]
    \centering
\begin{subfigure}{.48\textwidth}
  \centering   
\begin{tikzpicture}
 \node at (0,0) {\includegraphics[width=.8\linewidth]{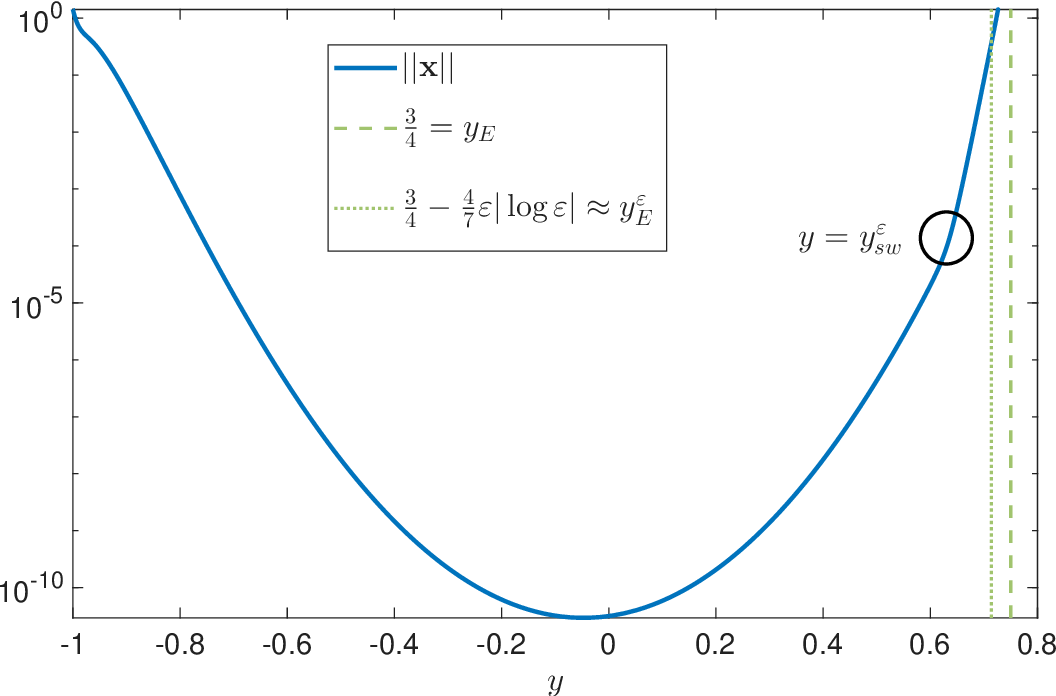}};
  \end{tikzpicture}
  \caption{}
\end{subfigure}
\begin{subfigure}{.48\textwidth}
  \centering
\begin{tikzpicture}
 \node at (0,0) {\includegraphics[width=.8\linewidth]{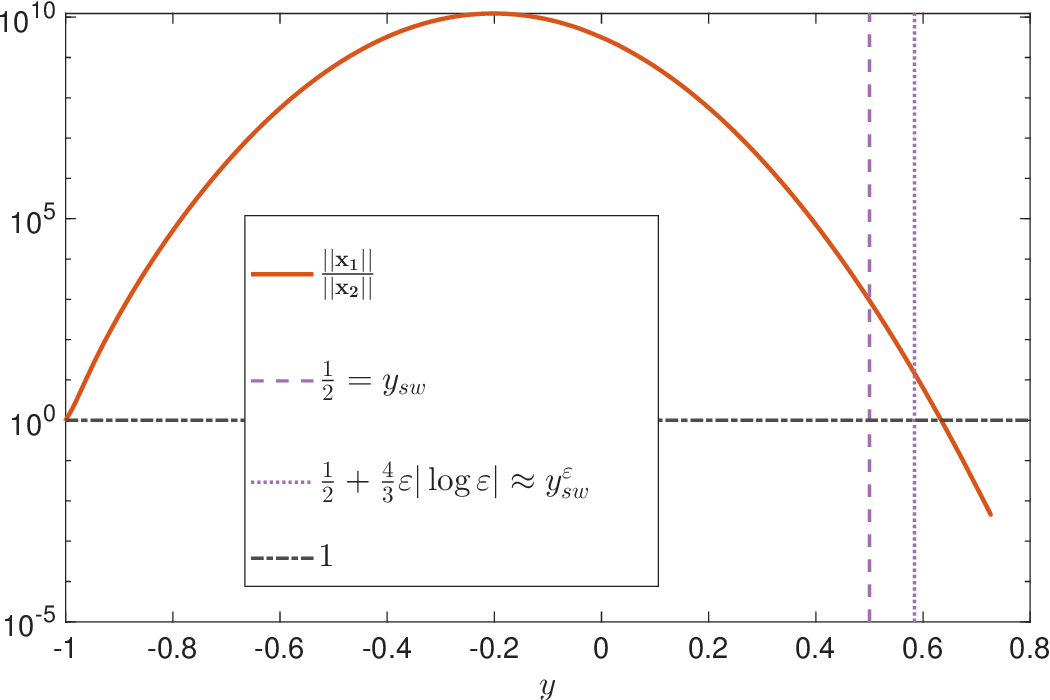}};
  \end{tikzpicture}
  \caption{}
\end{subfigure}
    \caption{Numerical simulation of \eqref{eq:example_2} for $\varepsilon=0.015$. (a) Evolution of $\norm{\mathbf{x}}$ during the entry-exit process, $y_E$ and an $\varepsilon$-approximation of $y_E^\varepsilon$ are indicated in green. As long as $y<y_{sw}^\varepsilon$, $\norm{\mathbf{x_1}}$ prevails, otherwise   $\norm{\mathbf{x_2}}$ prevails; notice indeed that $\norm{\mathbf{x}}$ starts to increase faster as $y=y_{sw}^\varepsilon$. (b) Evolution of $\norm{\mathbf{x_1}}/\norm{\mathbf{x_2}}$, $y_{sw}$ and an $\varepsilon$-approximation of $y_{sw}^\varepsilon$ are indicated in purple; $y_{sw}^\varepsilon$ corresponds to the intersection with the horizontal black line. Notice that the ratio of the norms is exponentially large in $\varepsilon$ away from where the norms are equal.\label{fig:numeric}}
\end{figure}

Figure \ref{fig:numeric} shows a numerical simulation of \eqref{eq:example_2}. These numerical results confirm our analytical results. In particular, the $\mathcal{O}(\varepsilon|\log \varepsilon|)$-corrections for $y_E$ and $y_{sw}$ allow us to derive good approximations for $y_E^\varepsilon$ and $y_{sw}^\varepsilon$, with $\mathcal{O}(\varepsilon)$-errors like in the usual settings. Moreover, even if there are no switching points close to $y_0$, $\norm{\mathbf{x}}$ shows a slightly different behavior close to this point because of the different mechanisms governing the initial dynamics. 

\section{Discussion}

The mechanism identified in this work is related to broader physical phenomena. For example, in hydrodynamic stability, eigenvalue-based modal analysis may fail to capture substantial transient amplification generated by the non-orthogonality of the underlying modes \cite{TreTreRedDri,Schm}. Another connection is with exceptional points in non-Hermitian systems, where the coalescence of eigenvectors gives rise to non-semisimple Jordan structures and polynomial contributions to the evolution of the orbits \cite{Hei}. In the slow-fast setting considered here, such polynomial effect, induced by the nilpotent component, acquires an additional interpretation: over the long $\mathcal{O}(\varepsilon^{-1})$ passage time, it produces the factor $\varepsilon^{-d}$ and ultimately shifts the exit point by $\mathcal{O}(\varepsilon|\log\varepsilon|)$. 

Finally, future work will investigate switching points in the degenerate case where $\lambda_i=\lambda_j$ at the spectral crossing. In this setting, the loss of transversality may lead to different asymptotic scaling laws and qualitatively new phenomena, including the emergence of multiple switching points.

\bigskip
\noindent\textbf{Acknowledgments.} ChatGPT (OpenAI, v. 5.6) was used to double-check and improve this work. The author takes full responsibility for the content of the manuscript. The author is a member and acknowledges the support of {\it Gruppo Nazionale di Fisica Matematica} (GNFM) of {\it Istituto Nazionale di Alta Matematica} (INdAM).

{\footnotesize
	\bibliographystyle{unsrt}
	\bibliography{biblio}
}

\end{document}